# The incompleteness of an incompleteness argument

Joachim Derichs

**Abstract**

Gödel's argument for the First Incompleteness Theorem is, structurally, a proof by contradiction. This paper intends to reframe the argument by, first, isolating an additional assumption the argument relies on, and then, second, arguing that the contradiction that emerges at the end should be redirected to refute this initial assumption rather than the completeness of number theory

# Section 1 – Introduction

It's the signature move of Gödel's argument to establish, based on the eponymous Gödel numbering, an ingenious mapping into number theory. But a mapping always has two sides, a preimage and an image, and we need to be as clear about the first as we are about the second. So what, exactly, is being mapped into number theory ?

On one side the mapping is already wonderfully clear, as it lands in first-order number theory. On the side of the preimage, much murkier. What we typically find there is a semi-formalism, assembled from different elements that are not well integrated: usually some recursion theory, often some naïve string set theory, always some highly uncommon notations to deal with encodings. All in all, a setup that seems dangerously complicated, and confusing.

This paper contends that we can do better, and wants to show it by restating the argument in a form that is more formally rigorous, and much more transparent.

In essence, we want there to be, on the preimage side as well, a theory phrased in the language of the predicate calculus. This will allow us to reconstruct the argument on a securely syntactic basis so that every step is realised as either a derivation inside a theory or a mapping between theories.

Among all the systems of notation currently in service, predicate logic is by far the best understood and most rigorously tested through long practical use. That we will be able to eliminate from the argument all elements not compatible with its rules brings several advantages:

> Maximum safety, as a standard in general use replaces niche formalisms like recursion theory and devices like corners or overlines that have no application outside metamathematics.

Maximum rigour, as a fully formal theory replaces an ad hoc assemblage of various notations at the origin of mapping, and concerning the target any suggestions of mapping to "numbers" rather than a theory are driven away.

Maximum transparency, as having the identical type of formalism on either side allows us to move from a complicated and idiosyncratic definition of the mapping to a much simpler one, in the form of a thoroughly normal and instantly readable homomorphism.

At the end of the process we will have shown the current messily mixed, residually informal setup can be recovered inside predicate logic, but not without changes. Although downstream many parts remain virtually the same, restating the antecedents does affect the substance of the result. Certain features that the traditional argument took for free are forced to assume their full cost. Under the strictest possible formalisation provided by unmodified predicate logic the features thus unearthed are recoverable only by making explicit an assumption that was implicit before in the structure of less reliable forms of notations. That the assumption is made in the course of Gödel-style arguments, either covertly or honestly, is not debatable; whether it is defensible can and will be debated, to an extent.

Once done clarifying, with a theory on both sides, we shall proceed to set up clean syntactic mappings between them:

**Definition 1.1:** A modelling transformation is a function m from the strings of a theory **T** to the strings of a theory **M** such that

1. m respects the roles of strings, mapping

    terms of **T** to terms of **M**,
    predicates of **T** to predicates (of the same arity) of **M**,
    sentences of **T** to sentences of **M**.

2. m commutes with logical connectors, quantification, and substitution, i.e.

| | |
|---|---|
| for all sentences φ | m(¬φ) = ¬m(φ) |
| for all sentences φ, ψ | m(φ ∨ ψ) = m(φ) ∨ m(ψ) |
| for all predicates P, (*x* free) | m(∃*x* P(*x*)) = ∃*x* m(P(*x*)) |
| for all predicates P, terms t | m(sub(P,t)) = sub(m(P),m(t)) |

3. m preserves provability, if **T** |- φ then **M** |- m(φ)

With an ironic nod to the semantic usage, we will call a consistent and complete **M** that receives a **T** under such a transformation a *model* for **T**.

The clarified preimage for Gödel's mapping should not only be a theory, but a theory of a peculiar kind. It needs to be literal about the strings of number theory in this technically precise sense:

**Definition 1.2:** A theory **M** is said to be meta to a theory **T** if all strings belonging to **T** are well-formed as primitive terms in **M**.

In a theory meta to number theory, the strings of number theory figure as closed terms. String predicates – e.g. the predicate of being a well-formed sentence – must first exist, unencoded, in a string theory before they can be projected, under a code, into number theory. The theory of number-theoretic strings is where string predicates live, in number theory they are only guests.

Interpretation, whether formal or informal, semantic or syntactic, is necessarily a binary relation, interpreting an A *as* a B, the image *as* the preimage. Before we can even think of interpreting arithmetic relations *as* string relations, we must secure possession of the string relations first.

So before Gödel's encoding argument can even begin, we need the following: At the origin of the mapping, a theory that is meta to number theory, and thus allows string predicates over the strings of number theory to be defined. The most important relation for the argument being provability, all the other relations only serve to support it.

Let **AR** stand for any first-order successor-based theory that contains Robinson's **Q**. **AR** could be number theory, or any suitably weaker subtheory.

What we need, then, for Gödel's argument to take off – or at least to take off in a clean, purely syntactic framework – is, minimally, this:

> **Definition 1.3: STRING** is a first-order theory that is meta to **AR** and contains a predicate PROV($x$) that directly defines provability for **AR**.

A distinct concept calls for distinct terminology. We say *directly* defines in order to emphasise again that the strings of the base theory appear in **STRING** nakedly, untouched by any encoding. Strings come just as they are. Literally literal.

For now, to maintain the symmetry between preimage and image, we shall limit the search for **STRING** to first-order theories. We promise to circle back to this limitation before the end.

The plan for the remainder of the paper is as follows:

(Section 2)  Examine **STRING**, the distilled preimage

(Section 3)  Examine concatenation, the crucial function in **STRING**

(Section 4)  Build modelling transformations between **STRING** and **AR**

(Section 5)  Recover a form of Gödel's Theorem

(Section 6)  Address the consequences of the restated Theorem

## Section 2 – Building **STRING**

This section examines the anatomy of **STRING**, the theory that is meant to constitute the preimage.

For the avoidance of misunderstanding, **STRING** is not and does not pretend to be what is often called 'the meta theory'. **STRING** is something much simpler: a concatenation-based, fully formal theory serving as the preimage for a mapping that the largely informal meta theory is used to describe. The language of **STRING** is first-order logic; the language of the background 'theory' is more or less lightly formalised English.

> **Definition 2.1:** Let *C* be the first-order language (with identity) consisting of a single binary function c(*x*,*y*), and for constants, strings over a given finite alphabet Σ.

Writing concatenation as a standard binary function is painfully inconvenient, which is why it is almost never done. Most presentations introduce specialised notations or conventions that incorporate the meta function of concatenation by juxtaposition, effectively putting "*xy*" in place of c(*x*,*y*). We will be avoiding these sleeker notations as a part of a principled commitment to staying within the notational means of predicate logic. In circumstances like ours, lack of polish can be a virtue. Clunky is good because it rattles whenever problematic moves are attempted.

For the following, note how among the predicate expressions from *C* there is a subset representing predicates that are primitive recursive in concatenation.

> **Proposition 2.2:** For **STRING** to exist it is sufficient that for any finite number of primitive recursive string predicates $P_i$ (i ≤ n) a theory **T** in *C* exists that
> (1) contains concatenation and string identity
> (2) decides the predicates $P_i$, so that for all strings s over Σ,
>     **T** |- $P_i$(s) or **T** |- ¬$P_i$(s)

Condition (1) requires the meaning of predicate expressions to be anchored by the c-function. What precisely is entailed by capturing concatenation will be the topic of the next section. For the duration of this section we will follow the example of all the textbook arguments for Gödel's Theorem and take native concatenation for granted. We will assume that the c-function behaves as expected, that it concatenates unencoded literals in just the way a naïve user would expect it to.

With the native c-function a given, the remaining content of the proposition covers well-trodden ground. We need only sketch a proof.

The intermediate step required is a predicate expression PROOF($x,y$) that formalises the string relation '*x is a proof sequence for the sentence y*'. As first shown by Gödel himself, there exists a predicate PROOF*($n,m$) in number theory that, as far as concatenation extends, could serve as a homomorphic image for PROOF. Moreover, PROOF*($n,m$) would be primitive recursive, and readily decidable. So we can simply argue that if PROOF is definable indirectly as an arithmetic relation, then it must be even more easily definable directly, without the extra complication of an encoding.

It is common knowledge that concepts like NUMERAL, VARIABLE, TERM, PREDICATE, SENTENCE are primitive recursive in concatenation; as is the relation of being an immediate consequence. Putting it all together, PROOF is quickly built.

So assuming only that the meaning of predicate expressions from *C* is anchored by concatenation, PROV($x$) := ∃$y$ PROOF($y,x$) will directly define provability.

To be clear, definability as used in the proposition is *not* intended to be definability according to this conventional definition:

**Definition 2.3**: A predicate P from an arithmetic theory **AR** defines_traditional a preformal string property **PROP**, if for all s∈Σ*, g an encoding function,

$$\models \textbf{PROP} \text{ is true for } s \Rightarrow \textbf{AR} \vdash P(g(s))$$

First of all, we do not want to encode. Instead, as already practised above, **STRING** means to formalise facts about naked literals. Hence our g is absent / the identical function. (To avoid confusion, we might choose to change the formatting for strings from the base theory, to **bold** or in some other uniform way, but that's it. No character replacement, and certainly no numeric encoding.)

Second, the way in which the conventional definition includes unexamined properties, along with "structures" that are semiformal at best, is profoundly unsatisfactory. It would be irresponsible to put any trust into unformalised ideas about string properties unless and until they have been refined into the predicates of an actual theory. The conventional definition has it backwards. One should not be referring theories to informal notions as the standard from which to take directions. We want theories to set the standard, and eliminate informal notions like **PROP**, as well as the accompanying truth talk.

Concatenation-based string theories represent much the clearest expression, the least slippery grip we have on the concept of a string property. String properties are so obviously at home there that it would in many ways be better to turn the definition talk around and say that the existence of a predicate in a string theory is what makes a well-defined string predicate. The predicates of string theories would then "define" by definition.

Before we had, for sentences φ of the arithmetic base theory,

$$\models \varphi \text{ is provable in } \textbf{AR} \Rightarrow \textbf{AR} \vdash \text{PROV*}(g(\varphi))$$

Now PROV *directly defines* provability, so that

$$\textbf{STRING} \vdash \text{PROV}(\varphi) \Rightarrow \textbf{AR} \vdash \varphi$$

Much simpler, no more g, and everything underpinned by thoroughly formal definitions. Since assuming c works at all, it works as well for the negated predicate, we also get a second line:

**STRING** $\vdash \neg PROV(\varphi) \Rightarrow$ **AR** $\nvdash \varphi$

Later on, when ready to start mapping to **AR**, we will be redefining representation as:

**STRING** $\vdash PROV(\varphi) \Rightarrow$ **AR** $\vdash PROV^*(g(\varphi))$

So eventually, we will get **AR** $\vdash PROV^*(g(\varphi))$ back, but will have eliminated all informal notions.

> **Definition 2.4**: A first-order theory is said to be strictly axiomatisable if is it axiomatised by a finite number of sentences or schemes such that any scheme can be summarised by a single sentence in higher-order logic.

The single sentence from higher-order logic provides a pattern against which axioms can be checked. Clearly, Peano arithmetic is strictly axiomatisable – its scheme is summarised in second-order logic. For naturally occurring theories, including of course **Q**, the definitions of 'strictly axiomatisable' and 'recursively axiomatisable' are coextensive. The adjusted definition is intended only to exclude artificial axiomatisations by pure enumeration or infinite indexing, which recursive axiomatisability might allow through.

> **Corollary 2.5**: For any consistent strictly axiomatisable theory **T**, if **STRING**$_T$ is complete, then
>
> **STRING**$_T$ $\vdash PROV_T(\varphi) \Leftrightarrow$ **T** $\vdash \varphi$

> **Definition 2.6**: A theory that contains concatenation is said to *directly express* any string relation that it decides.

**STRING**, if it were complete, would directly expresses provability for **AR**.

## Section 3 – On concatenation

In the previous section we have shown that the existence of a meta string theory that directly defines provability reduces to the existence of a theory that contains concatenation.

Speaking of concatenation, we have to be clear which form of concatenation we mean. There are two related but essentially different forms of concatenation. Both have their origin in the naïve concept of concatenation, but they do address different aspects.

The first form of concatenation is implicit. It is well-known and well-understood. Concatenation is formalised as an implicit binary function not unlike the Peano successor function. A typical sentence of such a theory would be the axiom of associativity: $c(x,c(y,z)) = c(c(x,y),z)$. An example of a weak theory of implicit concatenation, with most of the axioms going back to Tarski, is contained in Grzegorczyk [2005]. Stronger versions add an induction scheme, where the successor equivalent of a string is the string extended by one letter. The strongest version is defined, just as it is for arithmetic, by its second-order axiom of induction.

Implicit concatenation theories are characterised by the fact that they allow for only finitely many constants, one for each letter of their alphabet. Let's say that there are only two letters, a and b. The string "ab" would then not be a term of the theory. Implicit concatenation theories are therefore unable to make, let alone decide, statements of the form "c(a,b) = ab".

One could paraphrase this by saying that implicit concatenation theories express naïve concatenation only up to typographical permutation. Implicit concatenation only captures structural properties that are independent of an instantiation in concrete strings. It is not able to prove or refute that the concatenation of two given concrete strings is a third concrete string.

The second form of concatenation is much less well understood. We will call it extensional concatenation. It can be thought of as implicit concatenation made concrete, or fully interpreted.

People, and mathematicians especially, use it fluently, and daily; but very rarely do they stop to examine it.

Suppose we start again with only two letters, and want a language that will prove all and only the 'true' atomic concatenation statements (as naïvely understood) about these letters. For letters a and b that would mean that the language should prove e.g. c(a,b) = ab and c(a,ab) = aab, and disprove c(a,a) = a.

It should have become clear by now that extensional concatenation is a different and distinct concept from what is captured in theories of implicit concatenation. Extensional concatenation fixes the meaning of strings that implicit concatenation leaves to vary.

The point we have been building up to is that Gödel-style meta arguments absolutely require extensional concatenation. When the argument turns on a concrete string encoding another concrete string, representation up to typographical permutation is obviously not good enough.

Hence this first rough statement of an assumption:

> **Gödel's Assumption (rough) [3.1]:** Extensional concatenation is a well-defined function; and as such, primitive recursive.

The standard argument assumes primitive recursiveness for plenty of string relations. PROOF and all the other relations involved in its construction can only be primitive recursive if native concatenation, which lies at the base of their definition, is well-defined and primitive recursive, too.

There is no way of avoiding the use of extensional concatenation in Gödel-style arguments, as there is no way of talking about the strings of an underlying theory without it. Given the use, we would want to know that the naïve concept of extensional concatenation can be trusted. The prime means of demonstrating a clear understanding, and inspire trust, would be to present a consistent formalisation.

Now we want to avoid talk of recursion because recursive formalisations for native string properties are almost inevitably circular, assuming concatenation by juxtaposition as part of the definition.[1]

Due to reservations about recursion theory, and more important, to fit underneath **STRING**, we expect extensional concatenation to be formalisable in predicate logic instead:

> **Gödel's Assumption (still rough) [3.2]:** Extensional concatenation is formalisable. In other words, there exists a theory **CONCAT** containing a binary function c that correctly expresses extensional concatenation for atomic sentences over a finite alphabet $\Sigma$, so that for all $s,t,u \in \Sigma^*$,
>
> **CONCAT** $\vdash$ c(s,t) = u iff u is the concatenation of s and t

As a ubiquitous meta function, it is extremely easy for informal uses of extensional concatenation to intrude into an argument. The best protection against such intrusions is what we are already planning to do: Restate the preimage as a theory, and the entire interpretation process as a mapping from theory to theory. Reconstruct the argument on a purely syntactic basis, as an interpretation free of semantics. To that end, we will now try to formalise extensional concatenation. This will feel awkward as it is not something the predicate calculus is designed to do. The very awkwardness will show that the protections are working.

Gödel's Assumption, as stated, is infuriatingly vague. 'Correctly expresses' is just a pretentious way of saying 'behaves as one would naïvely expect concatenation to behave'. The phrase 'is the concatenation of' assumes command of the very thing we want do more about than just assume. But the task we face is tricky: How does one articulate the inarticulate ? How do we even know that there is something there to be expressed ?

---

[1] In recursion theory and related formalisms, definitions of native concatenation typically look like this: c($x,y$) = "f($xy$)", where f relies on concatenation by juxtaposition (Note the absence of a comma of the right-hand side). Rather than being openly declared, concatenation forms part of the formalism's structure, which renders its presence silent and near invisible.

In order to get a firmer grip let's try

> **Definition 3.3:** **CONCAT$_n$** is a consistent, complete and strictly axiomatisable first-order theory that directly defines extensional concatenation for strings up to length n over a finite alphabet Σ.

Writing out examples of **CONCAT$_n$** is straightforward if tedious. To illustrate, we present **CONCAT$_2$**, a theory consisting of six constants, a single predicate letter C, and the following axioms:

(As the assumption that concatenation should be a total function has become problematic, we switch to predicate-style notation to accommodate partial functions.)

Write $U_1(x) := x = \mathbf{0} \lor x = \mathbf{'}$

$U_2(x) := U_1(x) \lor x = \mathbf{00} \lor x = \mathbf{0'} \lor x = \mathbf{'0} \lor x = \mathbf{''}$

| | |
|---|---|
| (A1) $\forall x\ U_2(x)$ | There is an explicit, pre-declared finite domain. This in itself is nothing unusual, similar limitations are implicit in the workings of any real-world computer. |
| (A2) $\forall x,y\ (\neg U_1(x) \lor \neg U_1(y)) \to \neg \exists z\ C(x,y,z)$ | Concatenation in **CONCAT$_n$** is only a partial function, not defined beyond a given upper limit (again analogous to real-world computing). |
| (A3) $C(\mathbf{0},\mathbf{'},\mathbf{0'})$ | Within the strictures of the predicate calculus it is |
| (A4) $C(\mathbf{0},\mathbf{0},\mathbf{00})$ | surprisingly difficult to give correct extensional |
| (A5) $C(\mathbf{'},\mathbf{0},\mathbf{'0})$ | values in any other way than by listing all |
| (A6) $C(\mathbf{'},\mathbf{'},\mathbf{''})$ | instances individually. |

(A7) ¬(**0** = **'**)              Necessary to prevent trivialising models where letters are equated to each other. Collapsed concatenation with all letter set to be equal is (nearly) successor.

The theory **CONCAT**$_2$ is evidently decidable. The first axiom, by limiting the domain to a finite number of explicit strings, ensures that quantifiers can be eliminated in favour of atomic statements, in a finite number of steps. Axioms (A2) – (A6) together make sure that the list of atomic statements in C to be checked against is exhaustive, and in turn quickly decided. Atomic statements in = are resolvable into finite combinations of atomic statements in C and equalities over letters only. Last, atomic (in)equalities over letters are decided with the help of (A7).

It seems reasonable to grant the intelligibility of the idea of constructing theories **CONCAT**$_n$ on roughly this pattern even for n larger than the physical capacity of any human or computer. As categorical theories that directly define, the **CONCAT**$_n$ perfectly embody the concept they present. In their albeit limited realms, they *are* concatenation.

There is only one final complication: the self-reflexivity of extensional concatenation. Due to its circular nature, the construction of **CONCAT**$_2$ we gave does not, strictly speaking, prove the existence of languages **CONCAT**$_n$: The ability to construct **CONCAT**$_n$ presupposes string handling abilities tantamount to **CONCAT**$_m$, for some m > n. It therefore seems preferable to phrase the conclusion that these theories exist not as a proposition, but a postulate. Moreover, to avoid the implication that an infinite index i: n → **CONCAT**$_n$ exists, we are going to phrase the postulate as a scheme with parameter M:

      **Postulate Scheme [3.4]:** For n < M, **CONCAT**$_n$ exists.

Let M be a massively large finite number. For n < M, we concede that **CONCAT**$_n$ exists. The point of M is not make a specific number the limit, but to make it clear that there always is a limit. Every introduction of **CONCAT**$_n$ into a derivation comes with a side constraint that n < M. Let's call M the perimeter of concatenation. Inside the perimeter, everything works as (naïvely) expected. Most of the time, work goes on normally, and the perimeter can be ignored. Only when a contradiction arises is the perimeter constraint activated to absorb it. Sometimes, the argument can be repaired by shifting to a larger perimeter, M+X. So for example, a bounded diagonal argument that produces an element provably unlike any inside the perimeter can be repaired by shifting to a larger perimeter that includes the new element. At other times, the argument cannot be repaired. Diagonal arguments based on an infinite index cannot be restated to be compatible with a perimeter constraint. Such arguments can only be made by appealing to a more powerful postulate.

With the theories **CONCAT**$_n$ in mind we are now able to clarify

> **Definition [3.5]:** A first-order theory, styled **CONCAT**, is said to contain concatenation if it contains a function c that satisfies the following definition scheme with parameter M:
> For n < M, for all strings s, t, u in $\Sigma^n$
> **CONCAT** |- c(s,t) = u iff **CONCAT**$_n$ |- C(s,t,u)

We keep c as a function on the left-hand side in order to ensure continued compatibility, and interchangeability of predicates, with theories of implicit concatenation.

> **Gödel's Assumption (clarified) [3.6]: CONCAT** exists.

What began as something slippery and semantic – a claim about the legitimacy of using extensional concatenation in meta arguments – has now been turned into something purely syntactic, and thereby, unambiguous. Deciding atomic statements the same way as the theories **CONCAT**$_n$ is a pretty minimalist definition of what extensional concatenation should mean. Whatever one takes extensional concatenation to mean, however one tries to pin it down to a precise meaning, it is hard to see how any proposed formalisation of concatenation could be considered successful without meeting at least this lenient condition. As a formal statement of Gödel's Assumption it seems more than fair.

What is to stop us from disposing of Gödel's Assumption by simply proving it, delivering a theory that contains concatenation ? Can it really be so difficult to construct and axiomatise extensional concatenation ? Well, it's not so simple, it turns out.

Concatenation is an unusual function. To appreciate just how unusual, observe that the candidate that first comes to mind when one thinks about compacting the lists of explicit extensional values into a single axiom – "$\forall x,y \; c(x,y) = xy$" – is certainly not well-formed as a first-order formula, or indeed in any conventional logic. Given how it allows quantification to reach inside terms, inside the "$xy$", in a way that mashes meta and object level, it would not be surprising if the leading candidate for *Axiom of Concatenation*, once plausibly fleshed out into a non-standard theory, would turn out to be inconsistent.[2] Is this the notion that the average logician unconsciously applies in appealing to their intuition of concatenation ?

Without the trick of quantifying inside terms, **CONCAT** is in some trouble. The second idea that comes to mind – interpreting $xy$ as an (implicit) functional $x*y$, or even more explicitly as a function, $c*(x,y)$ – also fails. It leads to an immediate regress in the definition: $\forall x,y \; c(x,y) = c*(x,y)$ is not helpful. Thanks to the stubborn rigour of the predicate calculus, the informal function of concatenation by juxtaposition that we must inevitably make use of to operate any theory is not so easily hijacked for undercover work on the inside of an axiom.

So it seems that in order to meet the demands of expressing the concept of extensional concatenation, **CONCAT** would have to fall back on the unloved, third-best idea of listing explicit values – which is not without complications, either. Without the trick, or an implicit functional, what stands on one side of the equation defining concatenation must be a primitive term: the **0'** in c(**0**,') = **0'** must be primitive.

---

[2] More than a suspicion. Inconsistency for a non-standard formalism built to contain the *Axiom of Concatenation* is indeed provable. Unpublished paper.

(This is perhaps even easier to see in predicate-style notation, C(**0',0'**). When functionals have been eliminated altogether, there is no doubt that **0'** can only be a primitive term.) While manageable over finite domains, and perfectly workable for the definition of the theories **CONCAT**$_n$, when **CONCAT** tries to gather all the finite theories together, it has to start dealing with infinitely many primitive terms simultaneously. All of them must be handled, and fixed in their meanings by syntactic means only.

Formally, **0'** is a constant, a primitive term not analysable by its own language. A single symbol, without interior divisions. At the same time, the definition of constants as strings over $\Sigma$ presupposes a shadow theory that analyses them as non-primitive functional terms in concatenation. But we cannot forever rely on outsourcing work to the shadow theory, formalisation is not complete without knowing that we retain the capacity to internalise the work done by any actively involved third parties. The challenge for **CONCAT**, in a way, is to try and catch its own shadow. Put less poetically, in the mundane terms of computing: **0'**, as a single, undivided symbol, must have its own entry on the code table. This is a problem, as for all real-world computers the code table is finite.

Is the following operation computable ? Challenged by three objects, decide whether one is the concatenation of the other two. The hard part is formatting the problem, recognising the objects as strings. The inputs have to be parsed, someone or something has to break through the opacity of strings presented as a single symbol, and delineate letters. The rest, once parsing is complete, is stupendously easy.

At first, it may be difficult to see what the problem is. The hard part may not seem hard at all. Parsing is something that everyone who has learned to read an alphabetic language does involuntarily, unprompted. One hardly notices the effort, and finds it much more difficult to stop with the parsing than to parse. But stop is what we must do in order to be able to examine the operation, and really understand what is involved in extensional concatenation.

Once letters have been delineated, answering the challenge is as simple as deciding sentences in **CONCAT$_n$**. Yet parsing is not equivalent to deciding sentences in a given **CONCAT$_n$**, it is equivalent to setting up the **CONCAT$_n$**.

Deciding sentences inside **CONCAT$_n$** is trivial. Setting up any one theory **CONCAT$_n$** is still almost trivial. A general mechanism for setting up all the **CONCAT$_n$** is no longer trivial. It would, in fact, require infinitary powers.

The crux is to get from this

| | **0** | | **'** | | (on two consecutive squares)

to this

| | **0'** | | (on one square)

and vice versa.

Evidently, a machine able to effect concatenation and its inverse for all strings over a finite $\Sigma$ would have to be able to accept and process infinitely many distinct string objects. Translated into the language of Turing machines, **CONCAT** presupposes the ability of accepting infinitely many distinct symbols on a single square of tape. According to the standard definition of computability **CONCAT** is thereby disqualified as a computable function.

A theorem that takes in an uncomputable function, ready-made, is not going to have any difficulty proving uncomputability: Showing that functions derived from **CONCAT** are uncomputable, while true, fails to show anything about the computability of unrelated functions.

Though a blow, this is not a knockout yet because under an alternative definition of computability broader than Turing's, **CONCAT** might still be found to be a viable function. The decision whether **CONCAT** could exist as a non-standard function, provably diagonal to axiomatisations in standard predicate logic, is something we will have to return to.

## Section 4 – Mappings between **STRING** and **AR**

We are now able to put together **STRING**, the preimage: **CONCAT** is defined, and we already know that all string relations necessary for PROV are definable once **CONCAT** is given. Assuming only that **CONCAT** exists, then so does **STRING**.

For most intents and purposes, **CONCAT** already is **STRING**. The only bit missing is a means of deciding primitive recursive predicates.

There are various ways of axiomatising recursion in predicate logic, characterised by a trade-off between parsimony and elegance.

One can choose parsimony, by sticking closely to the concept of primitive recursion, with axioms that reproduce the mechanics of recursive enumeration in predicate logic. While certainly possible, this leads to cramped and unwieldy definitions.[3]

---

[3]Axiomatisation the parsimonious way proceeds by introducing < as an additional relation. The ordering relation needs to mesh with concatenation so that it contains the partial ordering by string inclusion, i.e. $c(x,y) = z \rightarrow ( x < z \land y < z )$, with axioms just strong enough to prove the finite definability of $x < s$, i.e.

For any string s over $\Sigma$,

**STRING** $\vdash x < s \leftrightarrow ( x = s_1 \lor x = s_2 \lor \ldots \lor x = s_m )$

where the $s_i$ are all strings over $\Sigma$ shorter than s

From finite definability it is then possible to derive equivalences of the form:

$P(x) \leftrightarrow A(x) \lor ( \exists x_i < x\ B(x,x_i) \land P(x_1) \land \ldots \land P(x_n) )$

where $A(x)$ is unary, and $B(x,x_i)$ is n+1-ary; both limited to bounded quantification

These provable equivalences are immediately recognisable as reproducing primitive recursion in the context of predicate logic.

Or one could opt for elegance, the transparency of intuitive definitions, but would then have to use axioms that are evidently stronger than direct equivalents of primitive recursion.[4]

Struggling for parsimony, however, goes unrewarded: All axiomatisations, no matter their degree of strength, incur the same debt by presupposing the existence of **CONCAT**.

For the progress of this paper the details of the axiomatisation by which **STRING** is given are therefore immaterial. The only fact that needs to be remembered is that the extension required to turn **CONCAT** into **STRING** is strictly axiomatisable.

Having established the preimage, and knowing the image to be a first-order arithmetic theory, we can move on to mappings between them. All the assumptions are in place for us to be able to connect **STRING** with **AR** by way of two lemmas.

> **Embedding Lemma for String Theories [4.1]:** For any strictly axiomatisable theory **T**, if **T** is consistent, and **STRING$_T$** is consistent and complete, then **STRING$_T$** contains a model of **T**.
>
> Proof: Let predicates P($x$) map to PROV(sub(P,$x$)), and terms identically. We show by induction on the composition of formulas that the mapping describes a modelling transformation. The only difficult case is **T** |- ¬φ ⇒ **STRING$_T$** |- ¬PROV(φ).
>
> > Assume **T** |- ¬φ. By consistency of **T**, **T** |/- φ.
> >
> > From completeness of **STRING**, ¬PROV directly expresses, so **STRING** |- ¬PROV(φ).
> >
> > The axiomatisability of **T** is required for the definability of PROV in **STRING$_T$**.

---

[4] By taking what induction has to offer one can enjoy the benefits of more fluent definitions for the building blocks of PROV:

(N1) NUMERAL(**0**) ("0 is a Numeral")

(N2) ∀$x$ ( NUMERAL($x$) → NUMERAL(c($x$,')) ) ("If $x$ is a Numeral, then $x$' is a Numeral")

(N3) ( P(**0**) ∧ ( ∀$x$ NUMERAL($x$) → ( P($x$) → P(c($x$,')) ) ) ) → ( ∀$x$ NUMERAL($x$) → P($x$) ) ("Nothing else is a Numeral")

A similar trio of clauses – consisting of a base clause, a construction clause, and an inductive closure clause – can axiomatise any recursively constructible string predicate. The clauses can either be introduced individually for each predicate, or more likely, derived from a general induction scheme.

**Gödel's Lemma [4.2]:** If **AR** is consistent, complete and strictly axiomatisable, **STRING** consistent, then **AR** contains a model of **STRING**.

Proof: Technically speaking, what needs to be shown is that when concatenation-based relations from **STRING** are associated with the right successor-based relations from **AR** by way of a Gödel numbering, then the resulting mapping is a modelling transformation.

This is a restatement of Gödel's groundwork for the First Incompleteness Theorem, accepting the substance with only minor changes: Instead of trying to map, on the preimage side, from informal properties or recursive theory or sets of strings, we are now mapping from the predicate strings of **STRING**. Instead of a largely unformalised meta-relation termed 'representation' we now have m(·), a provability-preserving homomorphism between first-order theories.

For all relations R of **STRING** it has to be shown that if **STRING** proves a sentence φ in R, then **AR** proves m(φ). Among the relations, we are really interested only in the first and last, c and PROV, but to reach the latter, all the relations in between have to be retraced.

First to be mapped is concatenation. As far as the c-function stretches, it can be shown that the targets envisaged in the standard argument are able to represent it.

Assume for simplicity that the target of c is a function, m(c) = c*. To even stake out the claim of interpretation for c*, we would need:

$$\textbf{CONCAT} \vdash c(s,t) = u \Rightarrow \textbf{AR} \vdash c^*(m(s),m(t)) = m(u)$$

while without the Assumption we have only:

For $n < M$, for all strings s, t, u in $\Sigma^n$

$$\textbf{CONCAT}_n \vdash c(s,t) = u \Rightarrow \textbf{AR} \vdash c^*(m(s),m(t)) = m(u)$$

Now it is true that far as finite claims of interpretation can be stated they can then also be shown to obtain. However, mapping atomic claims one by one is still not very promising, so on top of existing, we would want **CONCAT** to be axiomatisable. In order to get further than finitely far we would like axioms for c, so that we could then prove their translates in **AR**. Alas, we have none.

Failing axiomatisation, the other option is (second-order) meta-induction. For every letter a in Σ, right-concatenation $c_1 = c(x,a)$ and left-concatenation $c_2 = c(a,x)$ have to be shown to preserve provability under m, **STRING** $\vdash c_i(s) = u \Rightarrow$ **AR** $\vdash c_i*(m(s)) = m(u)$. Induction on the length of s and u would succeed, but does – obviously – require the assumption of **CONCAT**. So there is hole in the lemma, and it is filled with assumptions about **CONCAT**.

Next, all primitive recursive string relations required for PROOF, as far as the c-function will carry them, are shown to be homomorphic to their targets among arithmetic relations. With completeness, the homomorphism would extend to general relations like PROV.

As before, the axiomatisability of **AR** is required for the definability of PROV in **STRING**.

Let f and g be modelling transformations. Let g map **STRING** into **AR**, let f map **AR** into **STRING**, and for sentences φ put φ** = f(g(φ)). (Choice of letters is not accidental: Gödel's numbering g can be recovered as a mapping between terms induced by our g.)

We now have the ingredients ready for proving that the predicates of **STRING** under f o g have a fixed point. This follows generically from, for instance, category theory. We omit proofs for these familiar and easily verified conclusions.

> **Round Trip Lemma [4.3]**: STRING $\vdash \varphi \leftrightarrow \varphi$** (provided **AR** and **STRING** are both consistent and complete, and **AR** is strictly axiomatisable).

The cyclical image φ** of φ in **STRING** must be if not equal then at least equivalent to φ.

Section 5 – The reconstructed Theorem.

> **Definition 5.1:** The predicate TRUE is a truth predicate for a theory **STRING**$_T$ meta to **T**,

if for all sentences φ, **STRING**$_T$ |- TRUE(g(φ)) ↔ φ

where g is a modelling transformation from **STRING** into **T**.

What follows is a reconstruction or restatement of Gödel's Theorem within the framework we have developed so far. Though evidently not identical, it aims to stay as close as possible to the spirit of the original.

**Gödel's Theorem (according to Gödel) [5.2]:** If both **AR** and **STRING** are consistent and complete, and **AR** strictly axiomatisable, then **STRING** contains a truth predicate.

Proof: We show that PROV is a truth predicate.

**STRING** |- PROV(g(φ)) ⇒ **STRING** |- φ.

As PROV expresses provability for **AR**, **STRING** |- PROV(g(φ)) implies **AR** |- g(φ).

By embedding, **AR** |- g(φ) implies **STRING** |- f(g(φ)) = φ*.

Hence **STRING** |- φ (Round tripping).

**STRING** |- PROV(g(φ)) ⇐ **STRING** |- φ

**STRING** |- φ implies **AR** |- g(φ).

As PROV expresses provability for **AR**, **STRING** |- PROV(g(φ))

Predictably next,

**Tarski's Theorem [5.3]:** No consistent and complete theory can contain a truth predicate.

Proof: By applying to itself the predicate Fssb(P) ≡ ¬TRUE(g(sub(P,g(P)))), (F for false, ssb for self-substitution).

Gödel's cleverly self-referential sentence can be recovered as g(Fssb(Fssb)), the image in **AR** of the sentence Fssb(Fssb) that would break **STRING**.

**Corollary 5.4:** Either that

**STRING** exists, and is consistent and complete

or that

**AR** is consistent, complete, and strictly axiomatisable is false.

Proof: Taken together the two would prove the existence of a consistent and complete theory, **STRING**, with a truth predicate, contradicting Tarski's Theorem.

At first sight there appears to be an honest choice. The question seems to be on which of the initial assumptions should we pin the blame for a contradiction that emerges at the end of a long derivation. The theory **STRING** is a *prima facie* favourite for carrying the blame – it is an odd, untested and hugely ambitious theory. This conclusion becomes inevitable once we note that essentially the same argument can be made without involving arithmetic at all, simply by taking a second copy of **STRING** in place of **AR**. (It goes without saying that the argument can *not* be made with two copies of an arithmetic theory.)

**Gödel's Theorem (2$^{nd}$ approximation) [5.5]:** $STRING_{STRING}$, if consistent, complete, and strictly axiomatisable, contains a truth predicate.

Proof: Let $STRING_{STRING}$, a simple variation of $STRING_{AR}$, be a theory meta to a generic concatenation-based theory over a sufficiently large finite alphabet.

Take two copies of $STRING_{STRING}$, labelled $STRING_1$ and $STRING_2$. Each theory is, without any circularity, meta to the other.

If $STRING_1$ is axiomatisable, then by the same argument we used for **AR** a provability predicate PROV that applies to the underlying theory – $STRING_2$ is this case – can be defined. One way still takes the Embedding Lemma, the way back is much easier this time: the identical mapping. Round Trip then holds true, so the Theorem follows.

Although the investigation is not finished yet, there is one thing that we can already say with confidence: everything of importance in Gödel's argument happens on the preimage side. String theories are first in line to suffer incompleteness results. Only after they fall would these results also transfer to other kinds of theories. As the 2nd Approximation once again demonstrates, there can be Gödel-style incompleteness arguments without number theory, and any of its subtheories. But there can be no incompleteness arguments without **CONCAT**-based meta theories. At least on the preimage side, extensional string theories must always occur. From here on, as a consequence, arithmetic will be largely irrelevant to the investigation.

**STRING**$_{STRING}$ is in no way special. Any variation of **STRING** for different underlying theories consists of a version of **CONCAT** and varying definitions of meta predicates. The differing versions of **CONCAT** are contained in each other by a simple permutation of the alphabet; the meta predicates for one formal theory are definable by the same means as the predicates for any other. If we can successfully define a string theory meta to one axiomatisable theory, then we should be able to define a string theory meta to any other axiomatisable theory. Without loss of generality, we can therefore continue to speak of **STRING** without a subscript.

> **Corollary 5.6: STRING**, if it exists at all, and is consistent, cannot be both complete and strictly axiomatisable.
> Proof: By Tarski's Theorem, from the 2nd Approximation.

This is not quite the same as the trilemma that Gödel thought applied to number theory: From the three desirable qualities of consistency, completeness, and axiomatisability, pick any two. There is a forth option, that **CONCAT**, which was used to build **STRING**, does not exist. The contradiction that carries the proof is the same. But as a payoff for the cleanly syntactic presentation we have worked to achieve, we have an additional assumption for it to refute.

> **Corollary 5.7:** If **CONCAT** exists, then there exists a theory, **STRING**, that is essentially incomplete.

> Proof: Since the extension that expands **CONCAT** into **STRING** is strictly axiomatisable, the axiomatisability of **STRING** depends solely on the axiomatisability of **CONCAT.**

Except for the 'if **CONCAT** exists' condition, this finding of an essentially incomplete theory is the same as in traditional presentations. Incidentally, we could have reached an identical result by running the argument to this point with a concatenation theory on the image side, instead of **AR**.[5]

## Section 6 – Consequences

We are now in possession of a contradiction that refutes the conjunction of three assumptions, and we know what it would mean to deploy the contradiction to refute completeness / axiomatisability. What would it mean to deploy the contradiction to refute the new assumption on the board, Gödel's Assumption, the assumption of the right to use extensional concatenation in Gödel-style meta arguments, an assumption we have operationalised as the existence of **CONCAT** ? It would mean rejecting the existence of a theory that successfully envelops all the theories **CONCAT$_n$**. Approximating theories exist, but there would be, as it were, no limit to infinity.

In the end, there is a choice, though not an equal one. The situation is that there are two conflicting assumptions, neither of which appears to be provable outright. We have tried and failed to prove Gödel's Assumption from weaker or less problematic assumptions. On the opposite side we find the assumption of completeness for theories equipollent to number theory, e.g. number theory itself and comparably strong theories of implicit concatenation. Because for as long as we are also unable to prove the assumption of completeness – and it would be unrealistic to expect that to change – no hard

---

[5] This is a matter of some delicacy as it is known that adaptations of a Gödel-style incompleteness argument are available already for quite weak theories of implicit concatenation, e.g. Grzegorczyk [2005].

The point to note is that even though, confusingly, they bear the promising label 'concatenation', implicit theories of concatenation still have to be interpreted as extensional concatenation for any meta arguments to get underway. In this they are no different from number theory or **Q**, non-extensional theories that were interpreted as extensional concatenation for the sake of incompleteness arguments running in parallel. The results for implicit concatenation theories could be restated as we have restated the result for arithmetic, and would then also come back with an additional assumption about **CONCAT**.

contradiction will emerge. While the choice stands open, the completeness assumption, innocent though it may be, remains vulnerable to blame shifting.

Let's look first at the choice that makes less sense. Assuming a great deal of motivation it would remain just possible to continue insisting on the existence of **CONCAT**.

However, at this late stage it would effectively turn Gödel's Assumption into a postulate, Gödel's Axiom. The Axiom would assert the Assumption; or more generally, formalise assumptions that certain string predicate expressions successfully extensionalise into sets of strings. Through an axiom of this sort, naïvely taking descriptive expressions at face value, one could force into existence an infinitary string object, a theory in the extensional sense of a set of formal sentences.

The infinitary string object thus created would be custom-made to be unformalisable. A flat set of sentences, unaxiomatised and essentially unaxiomatisable. With the set pretending to contain the otherwise unreachable "truth" about extensional concatenation, we would be provably unable to summarise it by finite axioms, or a finite number of regular schemes.

Although possible, we should be clear about what taking this course of action would mean. Because of mutual representability by Gödel numberings and similar techniques it would mean that undecidability comes roaring back for a wide range of theories. Once (re)introduced as **CONCAT**, it would immediately spread to other theories, including some as weak as **Q**.

In the final analysis, to adopt Gödel's Assumption as an axiom would mean postulating essential incompleteness, and all the weirdness it entails, because one considers it to be a desirable state of the world.

Be that as it may, we stand at a fork, with two paths ahead. One is known to exhaustion, the other temptingly unexplored. For what remains of the paper, we will be taking the more interesting route, by treating the soft contradiction as a hard $\bot$, and then choosing to consider Gödel's Assumption refuted.

> **Gödel's Theorem (3rd approximation) [6.1]:** Extensional concatenation is not formalisable in first-order logic: For any theory **T** with strings over $\Sigma$ for constants, there exists an n and an atomic sentence $c(s,t) = u$ that is decided by **T** in a different way than it is by $\textbf{CONCAT}_n$.
>
> Proof: This is only a statement of the negation of Gödel's Assumption aka '**CONCAT** exists'.

Recall that, as a working assumption, **CONCAT** was defined to be a first-order theory. An obvious line of inquiry would be to ask whether higher orders of the predicate calculus, or other logics, would make any difference. It appears not. Readers are invited to verify that moving up to the second order does not materially change the outcome. The most promising route to an axiomatisation, the trick of quantifying inside terms, violates the standards of any conventional logic, second order as much as first. The presumption becomes that no consistent formalisation exists, in any conceivable type of formalism.[6]

To refute the Assumption is to conclude that there is no formalism that can contain concatenation, that decides all c-sentences the expected way. No formalism that can meet the low standard set by the definition of **CONCAT** for deserving to be called 'extensional concatenation'.

The new path begins by admitting that there is no infinite c-function, native string relations never were total, there is nothing useful on the side of the preimage for arithmetic relations to represent. With only uninteresting exceptions, interpretation fails for non-finite relations.

The new path continues to understanding that the absence of total extensions is not caused by theories failing to live up to some independent truth. It is caused by the fact that there is nothing to these notions of truth, nothing to formalise in the first place. Fragmented formalisations perfectly express that string reality itself is fragmented.

---

[6] The only formalisms that do appear to be able to capture extensional concatenation are non-standard and inconsistent (see footnote 2). This supports the conclusion that the idea of the function ranging infinitely wide, transcending all perimeters, is incoherent.

Although unformalisability will disappear, strange occurrences do not. Extensional string predicates will continue to provoke phenomena that can look superficially similar to the traditional weirdness created by essential incompleteness.

There are fragments that do not assemble into a whole. Fragmented predicates $P_n$, based on **CONCAT$_n$**, are admitted to exist, but these fragments do not assemble into total predicates P. Finite combinations of $P_n$ and $Q_m$ are possible, but, in general, no countably infinite unions.

>**Definition [6.2]:** For any P(*x*) from *C*, let $P_n$(*x*) be the expression that restricts P(*x*) to $\Sigma^n$, i.e. strings over $\Sigma$ of length $\leq$ n.

There are approximating sequences that converge over initial ranges, up to the edge of their perimeter, but fail to reach infinity. Specifically, there are diagonal predicates D where for finite N < M we can prove:

**CONCAT$_1$** |- $\forall x$ $D_1$(*x*)

**CONCAT$_2$** |- $\forall x$ $D_2$(*x*)

**CONCAT$_3$** |- $\forall x$ $D_3$(*x*)

...

**CONCAT$_N$** |- $\forall x$ $D_N$(*x*)

but we still cannot prove $\forall x$ D(*x*). No **CONCAT$_n$** can prove it for an unrestricted domain; nor would implicit theories of concatenation, where D(*x*) also forms part of the language. On the contrary, the expectation is that a plausible, sufficiently powerful theory of implicit concatenation would *reject* $\forall x$ D(*x*).

The inference

$$\text{"}\forall_i \mathbf{CONCAT}_i \vdash \forall x\, D_i(x) \Rightarrow \forall x\, D(x)\text{"}$$

is intuitively compelling, but false. Not so much because it is wrong – after all, all strings would eventually get covered by the enumeration –, but because without **CONCAT** we cannot even state it. With perimeters properly applied, all we can state is for arbitrarily large M

$$\forall_{i<M} \mathbf{CONCAT}_i \vdash \forall x\, D_i(x) \Rightarrow \forall x\, D(x),$$

which is not compelling at all.

One might paraphrase Gödel's assumption as the ability to quantify over perimeters, so that one can coherently speak of "$\forall$ M", or "$\forall_i \mathbf{CONCAT}_i$". This is equivalent to being allowed to enter another line into the enumeration above, after $\mathbf{CONCAT}_N \vdash \forall x\, D_N(x)$, reading '...'

While there still is weirdness, the important point is that it is localised. Unlike the traditional symptoms these strange tics are not infectious, they affect only extensional concatenation, and not any unrelated concepts, or theories. Number theory, in particular, is going to be totally immune.

We can enjoy a much more orderly world in which all well-defined problems are expected to be decidable, and belief in the existence of unformalisable mathematical truth is considered quaint.

**Gödel's Theorem (final) [6.3]:** Extensional concatenation is not a well-defined function.

In many ways, this conclusion is less a theorem than an exercise in boundary policing, an argument for which functions should be accepted as existing. By showing how bizarre the consequences of admitting extensional concatenation would be it urges to conclude that the idea has no merit. The better choice of Gödel's Axiom is to postulate not Gödel's Assumption, but its negation.

Gödel's incompleteness argument, though endlessly revealing about string theories, tells us very little about number theory. The original incompleteness result, the root cause of all the others, is the result for **CONCAT**. The contradiction that, when the existence of **CONCAT** is assumed, can be used to refute the completeness of number theory, and many other theories, is frankly imported from the preimage side. Incompleteness results are without exception the effect of projecting onto other theories the strange properties and failings of a hypothetically assumed **CONCAT**.

When we do conclude that **CONCAT** does not exist, the projections will instantly stop, and number theory, as well as all other theories able to model string predicates, will no longer be affected by any kind of unformalisability. (This includes also implicit theories of concatenation.) Once we reject Gödel's Assumption, no reason remains for believing in Gödel's original conclusion. Although absence of counterevidence does not constitute evidence, it is only natural to revert to completeness as the default assumption.

The lesson of Gödel's Theorem: Concatenation, used naïvely, can be just as treacherous as the membership relation $\in$. So in a way Gödel's argument does for concatenation-based theories what Russell's paradox does for set theory: For both types of theories, the naïve assumption that all predicate expressions ought to be able to have an extension turns out to be untenable despite its overwhelming intuitive appeal.